\documentclass[12pt, a4paper]{amsart} 

\usepackage{mathptmx} 
\usepackage{amsfonts,amsthm,amssymb,dsfont,stmaryrd,mathrsfs} 
\usepackage[leqno]{amsmath} 
\usepackage{array}
\usepackage{colortbl}
\usepackage{comment} 

\usepackage{graphicx,caption} 
\usepackage[pdfstartview=FitH]{hyperref} 
\usepackage{pdfpages} 
\usepackage{longtable,supertabular, multirow} 
\usepackage{diagbox} 
\usepackage{doi}
\usepackage{tikz, tikz-cd}
\usetikzlibrary{matrix,arrows,decorations.pathmorphing}

\usepackage[hmarginratio={1:1}, vmarginratio={1:1}, top=3cm, bottom=2cm, left=3cm, right=3cm]{geometry} 

\usepackage{setspace}
\allowdisplaybreaks[1]

\usepackage{mathtools} 

\theoremstyle{plain}
\newtheorem*{thm*}{Theorem} 
\newtheorem{thm}{Theorem} [section]
\newtheorem*{thm2*}{Theorem \ref{thmDecNilp}}
\newtheorem*{thm3*}{Theorem \ref{thmDecNilp2}}

\newtheorem{thmA}{Theorem}
\newtheorem{lem}[thm]{Lemma}
\newtheorem*{lem*}{Lemma}

\newtheorem{prop}[thm]{Proposition}
\newtheorem*{prop*}{Proposition}

\newtheorem{cor}[thm]{Corollary}
\newtheorem*{cor*}{Corollary}

\theoremstyle{definition}

\newtheorem*{defn*}{Definition}

\newtheorem*{conjecture*}{Conjecture}

\newtheorem*{exmp*}{Example}

\newtheorem*{prob*}{Problem}

\newtheorem*{ques*}{Question}
\newtheorem{remk}[thm]{Remark}
\newtheorem*{remk*}{Remark}

\def\nn{\mathbb{N}}

\def\qq{\mathbb{Q}}

\def\ff{\mathbb{F}}

\def\ker{{\rm ker}}

\def\GL{{\rm GL}}

\def\Tr{\operatorname{Tr}}

\def\mc#1{\mathcal{#1}}

\begin{document}

\title{Rational point counts arising in the modular isomorphism problem}



\author{Chih-Chiang Kao}
\address{Department of Mathematics, Institute of Science Tokyo, 2-12-1 Ookayama, Meguro-ku, 152-8550, Tokyo, Japan}
\email{kao.c.aa@m.titech.ac.jp}

\author{Wei-Liang Sun}
\address{Department of Mathematics, National Kaohsiung Normal University, Kaohsiung 82444, Taiwan}
\email{wlsun@mail.nknu.edu.tw}




\begin{abstract}
An explicit formula for the number of rational points on the common zero sets of certain multivariate polynomials over finite fields of characteristic $2$ is provided. As an application, finite $2$-groups of nilpotency class $2$ with cyclic center are determined by their modular group algebras. Consequently, the modular isomorphism problem has a positive answer for finite $p$-groups of nilpotency class $2$ with cyclic center for every prime $p$.
\end{abstract}

\maketitle


\section{Introduction}
\label{sec1}


Let $F$ be a field of characteristic $p > 0$. The modular isomorphism problem asks whether an isomorphism between the group algebras $F G$ and $F H$ of two finite $p$-groups $G$ and $H$ implies that $G$ and $H$ are isomorphic. It appeared in Richard Brauer's survey \cite{Brauer1963survey} on the general isomorphism problem. This holds for abelian groups, as proved in \cite{Deskins1956}. A counterexample was recently found in \cite{GarciaLucas-Margolis-delRio2022}, where the groups have the nilpotency class $3$. It is interesting to consider groups of nilpotency class $2$. Some positive partial results on the modular isomorphism problem have been obtained. The reader may refer to two nice serveys \cite{Sandling1985survey, Margolis2022survey}.


In this note, we are motivated by D. Garc\'{i}a-Lucas and L. Margolis' article \cite{GarciaLucas-Margolis2024} on the classification of $p$-groups of nilpotency class $2$ with cyclic center, since such groups can be seen as a starting point for studying $p$-groups of nilpotency class $2$. In their article, they proved that the modular isomorphism problem has a positive answer for odd $p$, and for $p=2$ when $x^2 + x + 1$ is irreducible over the field. For $p=2$, they also provide a criterion over finite fields in terms of bijective $F$-linear maps between certain sets of common zeros. We outline their approach below. 


Suppose that $G$ and $H$ are non-isomorphic finite $2$-groups of nilpotency class $2$ with cyclic center such that $F G$ and $F H$ are isomorphic group algebras for a field of characteristic $2$. A reduction result \cite[Corollary~B]{GarciaLucas-delRio2024} by Garc\'{i}a-Lucas and \'{A}. del R\'{i}o suggests that $F$ can be assumed to be a finite field $\ff_q$. According to~\cite[Lemma~3.10]{GarciaLucas-Margolis2024}, $G$ and $H$ must belong to two different types of groups (denoted by $Q$ and $R$ therein) in Y.K. Leong's classification in~\cite{Leong1979}. On the other hand, a certain quotient of ideals of $\ff_q G$ is determined by the set $V_{\ff_q}(f_r, h_r)$ of common zeros of two polynomials $$f_r = \sum_{i=1}^r x_i y_i \quad \text{and} \quad h_r = \sum_{i=1}^r x_i y_i (x_i + y_i)$$ in $\ff_q^{2r}$, while a certain quotient of ideals of $\ff_q H$ is determined by the set $V_{\ff_q}(g_r, h_r)$ of common zeros of two polynomials $$g_r = x_1^2 + y_1^2 + f_r$$ and $h_r$ in $\ff_q^{2r}$, where $r$ is a positive integer determined by both $G$ and $H$. By the proof of~\cite[Proposition~5.5]{GarciaLucas-Margolis2024}, an isomorphism between $\ff_q G$ and $\ff_q H$ induces a bijective $\ff_q$-linear map between $V_{\ff_q}(f_r, h_r)$ and $V_{\ff_q}(g_r, h_r)$. In their article, they showed that such a bijective map cannot exist if $r \leq 2$ by proving that $|V_{\ff_q}(f_r, h_r)|$ is greater than $|V_{\ff_q}(g_r, h_r)|$. This shows that $\ff_q G$ and $\ff_q H$ must be non-isomorphic group algebras when $r \leq 2$.

In this article, we will show that $|V_{\ff_q}(f_r, h_r)|$ and $|V_{\ff_q}(g_r, h_r)|$ are distinct for every $r$ by providing explicit formulas for both cardinalities. 

\begin{thmA}
\label{thmA: formula of size of rational points}
Let $q = 2^n$ for $n \in \nn$ and let $\ff_q$ be a finite field of $q$ elements. For $r \in \nn$ and $k \in \{0, 1, \ldots, r\}$, let $$g_{r,k} = \sum_{i=1}^k (x_i^2 + y_i^2) + \sum_{i=1}^r x_i y_i \quad \text{and} \quad h_r = \sum_{i=1}^r x_i y_i (x_i + y_i)$$ be two polynomials over $\ff_q$, where the sum of squares is omitted when $k=0$, and let $V_{\ff_q}(g_{r,k}, h_r)$ be the number of common zeros of $g_{r,k}$ and $h_r$ in $\ff_q^{2r}$. Then $$|V_{\ff_q}(g_{r,k}, h_r)| = \frac{q^r}{q^2} \left( q^r + (-1)^{nk} \left( q - 1 + \frac{\mc{I}_q(3)}{q} + 2^r (-1)^k \mc{I}_q(2) + \frac{4^r \binom{q}{3}}{q} \right) \right),$$ where $\mc{I}_q(d)$ is the number of monic irreducible polynomials of degree $d$ over $\ff_q$.
\end{thmA}

With the notation above, one has $g_{r,0} = f_r$ and $g_{r,1} = g_r$. It is not difficult to see from the above formula that $|V_{\ff_q}(f_r, h_r)|$ and $|V_{\ff_q}(g_r, h_r)|$ are distinct. By Garc\'{i}a-Lucas and del R\'{i}o's reduction theorem and Garc\'{i}a-Lucas and Margolis' criterion, we obtain the following conclusion.

\begin{thmA}
\label{thmA: MIP to 2-groups of nil class 2 with cyclic center}
Let $F$ be a field of characteristic $2$ and $G$ a finite $2$-group of nilpotency class $2$ with cyclic center. If $H$ is a group such that the modular group algebras $F G$ and $F H$ are isomorphic, then the groups $G$ and $H$ are isomorphic.
\end{thmA}

Together with \cite[Theorem~1.1]{GarciaLucas-Margolis2024} for finite $p$-groups with odd $p$, this yields the following result.

\begin{thmA}
\label{thmA: MIP to p-groups of nil class 2 with cyclic center}
Finite $p$-groups of nilpotency class $2$ with cyclic center are determined up to isomorphism by their modular group algebras.
\end{thmA}

In this article, the required techniques from algebraic number theory are presented in Section~\ref{sec2}. Section~\ref{sec3} provides the computational consequences that will be utilized in Section~\ref{sec4} to prove Theorem~\ref{thmA: formula of size of rational points}. 

\section{Preliminaries}
\label{sec2}

Throughout this article, $\ff_q$ denotes the finite field with $q = p^n$ elements for a prime $p$ and a positive integer $n$. We recall the following standard terminology from \cite[Section~10.3]{IreRos}. The {\it trace} from $\ff_{q}$ to $\ff_p$ is defined as the sum of the images of an element under the Galois automorphisms generated by the Frobenius automorphism: $$\Tr_{\ff_{q} / \ff_p}(x) = x + x^{p} + x^{p^2} + \cdots + x^{p^{n-1}}$$ for $x \in \ff_{q}$. If there is no confusion, we write $\Tr_{\ff_{q} / \ff_p}(x) = \Tr(x)$ for simplicity. We denote the canonical additive character of $\ff_q$ by $\psi$ which is given by $$\psi(x) = \zeta_p^{\Tr(x)},$$ where $\zeta_p$ is a primitive $p$-th root of $1$ in $\overline{\qq}$. Thus, $$\psi(x + y) = \psi(x) \psi(y)$$ for $x, y \in \ff_q$. Since $\psi$ is a nontrivial additive character, it follows that \begin{align} \sum_{x \in \ff_{q}} \psi(x) = 0. \label{eq: sum of images of psi} \end{align} 

Note that $\Tr(x) = \Tr(x^p)$ for all $x\in\ff_q$, and hence $\psi(x) = \psi(x^p)$. Moreover, because $x\mapsto x^p$ is the Frobenius automorphism of $\ff_q$, every element in $\ff_q$ has a unique $p$-th root. One has \begin{align} \psi(x) = \psi(\sqrt[p]{x}). \label{eq: p-th root under psi}\end{align} Moreover, the map $x \mapsto x^p - x $ is an $\ff_p$-linear map on $\ff_q$ with kernel $\ff_p$. So the image has size $q/p$. Since $x^p - x \in \ker(\Tr)$ and the size of $\ker(\Tr)$ is also $q/p$, one deduces that \begin{align} \ker(\Tr) = \{x^p - x \mid x \in \ff_q\}. \label{eq: kernel of Tr} \end{align}

A solution to a multivariate polynomial with coefficients in $\ff_q$ is said to be {\it $\ff_q$-rational} if its coordinates lie in $\ff_q$.


\section{Character Sums and Trace Conditions}
\label{sec3}

Throughout this section, $\ff_q$ is a field of characteristic $2$. The following two sums will be widely used in the rest of the discussion: for $a, b \in \ff_q$, \begin{align} S(a, b) & = \sum_{x, y \in \ff_q} \psi(a xy + b xy (x+y)) \quad \text{and} \\ T(a, b) & = \sum_{x,y\in\ff_q}\psi(a(x^2+y^2+xy)+b xy (x+y)). \end{align}

\begin{lem}
\label{lem: S(a,b)}
    For any $a,b\in\ff_q$, we have
    $$S(a,b) = \begin{cases}
        q^2&\text{ if }a = b = 0,\\
        q&\text{ if }a\neq 0, b = 0,\\
        q(1+N(a,b))&\text{ if }b\neq 0,
    \end{cases}$$
    where $N(a,b)$ is the number of distinct roots of $z^{3} + a z + b$ in $\ff_q$.
\end{lem}

\begin{proof}
    By the definition of $S(a,b)$, it is immediate that
    $S(0,0) = q^2.$
    Now, assume that $a\neq 0$ and $b = 0$. Then $S(a,0) = \sum_{x,y\in\ff_q}\psi(axy)$. If $x = 0$, then $\sum_{y\in\ff_q}\psi(axy) = q$. On the other hand, if $x\neq 0$, then $y\mapsto axy$ is a bijection on $\ff_q$. Thus, \eqref{eq: sum of images of psi} implies 
    $\sum_{y\in\ff_q}\psi(axy) = \sum_{y\in \ff_q}\psi(y) = 0.$
    Hence, $S(a,0) = q$.

    Finally, suppose $b\neq 0$. According to \eqref{eq: p-th root under psi}, we have 
    $$\psi(bxy^2) = \psi((\sqrt{bx}y)^2) = \psi(\sqrt{bx}y).$$
    So $S(a,b) = \sum_{x,y\in\ff_q}\psi((ax+bx^2+\sqrt{bx})y)$. For a fixed $x\in\ff_q$, it follows from \eqref{eq: sum of images of psi} that 
    $$\sum_{y\in\ff_q}\psi((ax+bx^2+\sqrt{bx})y) = \begin{cases}
        q&\text{ if }ax+bx^2+\sqrt{bx} = 0,\\
        0&\text{ if }ax+bx^2+\sqrt{bx} \neq 0.\\
    \end{cases}$$
    Hence
    $$\sum_{x,y\in\ff_q}\psi((ax+bx^2+\sqrt{bx})y) = q\cdot|\{x\in\ff_q:ax+bx^2+\sqrt{bx} = 0\}|.$$

    Let $z = \sqrt{bx}$. Then $ax+bx^2+\sqrt{bx} = \frac{1}{b} (z^{4} + a z^2 + bz)$. Hence $ax+bx^2+\sqrt{bx} = 0$ if and only if $z = 0$ or $z^{3} + a z + b = 0$. Hence, we have $S(a,b) = q(1+N(a,b))$.
\end{proof}


\begin{lem}
\label{lem: T(a,b)}
    For any $a,b\in\ff_q$, we have
    $$T(a,b) = \begin{cases}
        q^2&\text{ if }a = b = 0,\\
        q\cdot\psi(1)&\text{ if }a\neq 0, b = 0,\\
        S(a,b)\cdot \psi(a^3/b^2)&\text{ if }b\neq 0.
    \end{cases}$$
\end{lem}

\begin{proof}
    By the definition of $T(a,b)$, we have
    $T(0,0) = q^2.$
    Now, assume that $a\neq 0$ and $b = 0$. By~\eqref{eq: p-th root under psi}, $\psi(ax^2) = \psi(\sqrt{a}x)$ and $\psi(ay^2) =  \psi(\sqrt{a}y)$ for all $x, y \in \ff_q$. One has
    \begin{align*}
        T(a,0) = \sum_{x,y\in\ff_q}\psi(\sqrt{a}x+\sqrt{a}y+axy) = \sum_{x\in\ff_q}\psi(\sqrt{a}x) \sum_{y\in\ff_q}\psi((\sqrt{a}+ax)y).
    \end{align*}
    Fix $x\in\ff_q$. By~\eqref{eq: sum of images of psi},
    $$\sum_{y\in\ff_q}\psi((\sqrt{a}+ax)y) = \begin{cases}
        q&\text{ if }\sqrt{a}+ax = 0,\\
        0&\text{ if }\sqrt{a}+ax\neq 0.
    \end{cases}$$
    If $\sqrt{a}+ax =  0$, then $x = \sqrt{a}/a$, which implies that $\psi(\sqrt{a}x) =\psi(1)$. Thus, $T(a,0) = q\cdot\psi(1)$.

    Finally, suppose $b\neq 0$. By~\eqref{eq: p-th root under psi}, $\psi(a y^2) = \psi(\sqrt{a}y)$ and $\psi(b x y^2) = \psi(\sqrt{bx} y)$. Thus,
    \begin{align*}
        T(a,b) & = \sum_{x,y\in\ff_q}\psi(ax^2+\sqrt{a}y+axy+bx^2y+\sqrt{bx}y)\\
        & = \sum_{x\in\ff_q}\psi(ax^2)\sum_{y\in\ff_q}\psi((\sqrt{a}+\sqrt{bx}+ax+bx^2)y).
    \end{align*}
    By~\eqref{eq: sum of images of psi},
    $$\sum_{y\in\ff_q}\psi((\sqrt{a}+\sqrt{bx}+ax+bx^2)y) = \begin{cases}
        q&\text{ if }\sqrt{a}+\sqrt{bx}+ax+bx^2 = 0,\\
        0&\text{ if }\sqrt{a}+\sqrt{bx}+ax+bx^2\neq 0.
    \end{cases}$$
    It follows that 
    $$T(a,b) = q\cdot\sum_{x\in\Omega(a,b)}\psi(ax^2),$$
    where $$\Omega(a,b) = \{x\in\ff_q:\sqrt{a}+\sqrt{bx}+ax+bx^2 = 0\}.$$ 

    We next compute $\psi(ax^2)$ for $x\in\Omega(a,b)$. Let $z = \sqrt{a}+\sqrt{bx}$. Then $x = \frac{z^2-a}{b}$. Moreover, since $x\in\Omega(a,b)$, we have $z^4 = az^2+bz$. Hence,
    $$ax^2 = \frac{a}{b^2}(z^4-a^2) = \frac{a}{b^2}(az^2+bz+a^2) = \left(\frac{az}{b}\right)^2+\frac{az}{b}+\frac{a^3}{b^2}.$$
    Since
    $$\Tr(u+u^2) = \Tr(u)+\Tr(u^2) = 2\Tr(u) = 0,$$
    it follows that $\Tr(ax^2) = \Tr(a^3/b^2)$ and so $\psi(ax^2) = \psi(a^3/b^2)$. Moreover, $\sqrt{a}+\sqrt{bx}+ax+bx^2 = 0$ if and only if $z = 0$ or $z^3+az+b = 0$. Since $x\mapsto z = \sqrt{bx}+\sqrt{a}$ is a bijection on $\ff_q$, and $0$ is not a root of $z^3+az+b$, we have $|\Omega(a,b)| = 1+N(a,b)$, where $N(a,b)$ is defined in Lemma~\ref{lem: S(a,b)}. Therefore, $q\cdot|\Omega(a,b)| = q(1+N(a,b)) = S(a,b)$, which implies that
    $$T(a,b) = q\cdot\sum_{x\in\Omega(a,b)}\psi(a^3/b^2) = q\cdot|\Omega(a,b)|\cdot\psi(a^3/b^2) = S(a,b)\cdot\psi(a^3/b^2).$$
\end{proof}

We are going to determine whether $\Tr(a^3/b^2)$ is $0$ or $1$ for each pair $(a,b)$ with $b \neq 0$, and then determine $\psi(a^3/b^2)$. First, we need the following result on quadratic polynomials.

\begin{lem}
\label{lem: psi(d/c^2)}
    Let $f(z) = z^2+cz+d$ be a quadratic polynomial in $\ff_q[z]$ with $c\neq 0$. Then $f(z)$ splits over $\ff_q$ if and only if $\Tr(d/c^2) = 0$.
\end{lem}

\begin{proof}
    Let $z = cT$. Then $f(z) = (cT)^2+c(cT)+d = c^2(T^2+T)+d$. Since $c\neq 0$, $f(z)$ splits over $\ff_q$ if and only if $T^2+T+d/c^2$ splits over $\ff_q$, which is equivalent to saying that $T^2+T+d/c^2$ has a root in $\ff_q$. If $t \in \ff_q$ is a root of $T^2+T+d/c^2$, then $$\Tr(d / c^2) = \Tr(t^2 + t) = \Tr(t^2) + \Tr(t) = 2 \Tr(t) = 0.$$ Conversely, if $\Tr(d / c^2) = 0$, then $d / c^2 \in \ker(\Tr)$. According to~\eqref{eq: kernel of Tr}, $d / c^2 = t^2 - t = t^2 + t$ for some $t \in \ff_q$, as desired. 
\end{proof}

\begin{lem}
\label{lem: number of (a,b) with Tr(a^3/b^2) = 1}
    There are $q(q-1)/2$ pairs $(a,b)\in\ff_q\times\ff_q^{\times}$ such that $\Tr(a^3/b^2) = 1$.
\end{lem}

\begin{proof}
    Since $b\mapsto b^{-2}$ is a bijection on $\ff_q^{\times}$, we may set $c = 1/b^2$, and it suffices to count the pairs $(a,c)\in\ff_q\times\ff_q^{\times}$ such that $\Tr(a^3c) = 1$. Recall that there are $q/2$ elements of $\ff_q$ with trace $1$ and $\Tr(0) = 0$. So there are $q/2$ elements of $\ff_q^{\times}$ having trace 1.

    If $a = 0$, then $\Tr(a^3c) = 0$ for all $c\in\ff_q^{\times}$. If $a\neq 0$, then $c\mapsto a^3c$ is a bijection on $\ff_q^{\times}$. Hence, for each fixed $a\neq 0$, there are $q/2$ elements $c \in \ff_q^{\times}$ such that $\Tr(a^3c) = 1$, as desired. 
\end{proof}

The value of $\Tr(a^3 / b^2)$ can be determined by the number of $\ff_q$-rational roots of the cubic polynomial $z^3 + a z + b$. Before we prove this, we count the number of pairs $(a, b)$ for which such a cubic polynomial has a given number of roots.

Let $\mc{I}_q(d)$ denote the number of monic irreducible polynomials of degree $d$ over $\ff_q$. 
It is well known (see \cite[Corollary to Proposition~2.1]{Rosen2002}) that $$\mc{I}_q(d) = \frac{1}{d} \sum_{a \mid d} \mu(a) q^{\frac{d}{a}}$$ where $\mu$ is the M\"{o}bius function. Thus, $$\mc{I}_q(2) = \frac{1}{2} (q^2 - q) \quad \text{and} \quad \mc{I}_q(3) = \frac{1}{3} (q^3 - q) .$$ 

\begin{lem}
\label{lem: mc{N}_k}
    Let $\mc{N}_k$ be the number of pairs $(a,b)\in\ff_q\times\ff_q^{\times}$ such that $z^3+az+b$ has exactly $k$ distinct roots in $\ff_q$. Then
    $$\mc{N}_0 = \mc{I}_q(3)/q, \quad \mc{N}_1 = \mc{I}_q(2), \quad \mc{N}_2 = 0, \quad \mc{N}_3 = \binom{q}{3} / q.$$
\end{lem}

\begin{proof}
Let $x_1, x_2, x_3$ be the roots of $z^3 + a z + b$ in an algebraic extension of $\ff_q$. Since $b\neq 0$, all roots are simple, and hence $\mc{N}_2 = 0$. If $x_1$ is the only root in $\ff_q$, then \begin{align} \frac{z^3 + a z + b}{z-x_1} = z^2 + x_1 z + x_1^2 + a \label{eq: proof of lem: mc{N}_k} \end{align} is a monic irreducible quadratic polynomial over $\ff_q$ with $x_1$ as the coefficient of the linear term. Conversely, each monic irreducible quadratic polynomial uniquely determines such a cubic polynomial. Thus, $\mc{N}_1 = \mc{I}_q(2)$.   

Next, $\mc{N}_0$ is the number of monic irreducible cubic polynomials of the form $z^3 + a z + b$. Replacing $z$ with $z-\frac{c}{3}$ in $z^3 + c z^2 + d z + e$ yields $z^3 + a z + b$, where $a = d - \frac{c^2}{3}$ and $b = e + \frac{2c^3}{27} - \frac{cd}{3}$. This transformation preserves irreducibility. For fixed $a$ and $b$, the coefficients $d$ and $e$ are uniquely determined by $c$, which can be chosen arbitrarily from $\ff_q$. Hence, this correspondence is $q$-to-one, giving $\mc{N}_0 = \mc{I}_q(3)/q$.

Finally, the preceding $q$-to-one correspondence yields $\mc{N}_3 = \binom{q}{3}/q$, since $\binom{q}{3}$ counts the monic cubic polynomials with three distinct roots in $\ff_q$.
\end{proof}

Recall that $N(a,b)$ is the number of distinct $\ff_q$-rational roots of $z^3 + a z + b$, as defined in Lemma~\ref{lem: S(a,b)}. The following result characterizes the values of $\Tr(a^3 / b^2)$ in terms of $N(a,b)$, which cannot be $2$ by Lemma~\ref{lem: mc{N}_k}.

\begin{prop}
\label{prop: determine Tr(a^3/b^2)}
    Let $q = 2^n$ for $n \in \nn$, and let $(a,b)\in\ff_q\times\ff_q^{\times}$. Then
    $$\Tr(a^3/b^2) = \begin{cases}
        n &\text{ if }N(a,b) = 0 \text{ or }3,\\
        n+1 &\text{ if }N(a,b) = 1.
    \end{cases}$$
\end{prop}

\begin{proof}
    First, suppose $z^3+az+b$ has a root $x_1\in\ff_q$, which is nonzero since $b \neq 0$. Then $N(a,b) = 1$ or $3$. According to~\eqref{eq: proof of lem: mc{N}_k} in the proof of Lemma~\ref{lem: mc{N}_k} and Lemma~\ref{lem: psi(d/c^2)}, $N(a,b) = 1$ if and only if $\Tr((x_1^2+a) / x_1^2) = 1$, and $N(a,b) = 3$ if and only if $\Tr((x_1^2+a) / x_1^2) = 0$.

    We denote $a/x_1^2$ by $\omega$. Note that $$\Tr((x_1^2+a) / x_1^2) = n + \Tr(\omega).$$ By~\eqref{eq: proof of lem: mc{N}_k} again, $b = x_1^3 + a x_1$.  Since $$\omega^3 = (1 + \omega)^3 + (1 + \omega)^2 + (1 + \omega) + 1,$$ it follows that $$\frac{a^3}{b^2} = \frac{a^3}{x_1^6+a^2x_1^2} = \frac{\omega^3}{1+\omega^2} = \frac{\omega^3}{(1+\omega)^2} = (1 + \omega) + 1 + \frac{1}{1 + \omega} + \frac{1}{(1 + \omega)^2}.$$ We obtain from~\eqref{eq: kernel of Tr} that $$\Tr(a^3 / b^2) = \Tr(\omega) = n + \Tr((x_1^2+a) / x_1^2) = \left\{ \begin{array}{ll} n + 1 & \text{if } N(a,b) = 1, \\ n & \text{if } N(a,b) = 3. \end{array} \right.$$
    
    According to Lemma~\ref{lem: mc{N}_k}, $\mc{N}_1 = \mc{N}_0 + \mc{N}_3 = q (q-1) / 2$, which is the same as the number of pairs $(a,b)$ with $\Tr(a^3/b^2) = 1$ by Lemma~\ref{lem: number of (a,b) with Tr(a^3/b^2) = 1}. 
    If $n$ is even, the pairs with $N(a,b)=1$ account for all pairs with $\Tr(a^3 / b^2) = 1$. If $n$ is odd, the pairs with $N(a,b) = 0$ or $3$ account for all such pairs. In either case, $\Tr(a^3 / b^2) = n$ when $N(a,b) = 0$.
\end{proof}


\section{Rational Point Counts}
\label{sec4}

Let $r\in \nn$, and define the polynomials
$$g_{r,k} = \sum_{i=1}^k (x_i^2 + y_i^2) + \sum_{i = 1}^r x_i y_i, \quad h_r = \sum_{i = 1}^r x_iy_i(x_i+y_i),$$ where $0 \leq k \leq r$, with the convention that the sum of squares is omitted when $k=0$.
For a field $F$ and a subset $S\subseteq F[x_1,y_1,\ldots,x_r,y_r]$, we define 
$$V_F(S) = \{P\in F^{2r}:f(P) = 0\text{ for all }f\in S\}.$$

\begin{thm}
\label{thm: sizes of rationals points}
    Let $q = 2^n$, $r \in \nn$ and $k \in \{0, 1, \ldots, r\}$. 
    The number of $\ff_q$-rational points of $V_{\ff_q}(g_{r,k},h_r)$ is $$\frac{q^r}{q^2} \left( q^r + (-1)^{nk} \left( q - 1 + \frac{\mc{I}_q(3)}{q} + 2^r (-1)^k \mc{I}_q(2) + \frac{4^r \binom{q}{3}}{q} \right) \right),$$ where $\mc{I}_q(d)$ is the number of monic irreducible polynomials of degree $d$ over $\ff_q$.
\end{thm}

\begin{proof}
    For a function $u:\ff_q^{2r}\rightarrow\ff_q$, the characteristic function $\mathds{1}_{u = 0}:\ff_q^{2r}\rightarrow\{0,1\}$ of the zero set of $u$ is defined by 
    $$ \mathds{1}_{u = 0}(P) = \begin{cases}
        1&\text{ if }u(P) = 0,\\
        0&\text{ if }u(P)\neq 0
    \end{cases}$$
    for $P \in \ff_q^{2r}$.
    It follows from~\eqref{eq: sum of images of psi} that, for any $P\in\ff_q^{2r}$
    $$ \mathds{1}_{u = 0}(P) = \frac{1}{q}\sum_{a\in \ff_q}\psi(au(P)).$$ For $a,b \in \ff_q$, we simply denote $$s_i(a,b) = a x_i y_i + b x_i y_i (x_i + y_i)$$ and $$t_i(a,b) = a (x_i^2 + y_i^2 + x_i y_i) + b x_i y_i (x_i + y_i)$$ for variables $x_i , y_i \in \ff_q$.
    Hence,
    \begin{align*}
        & |V_{\ff_q}(g_{r,k},h_r)| \\ = {} &  \sum_{P\in\ff_q^{2r}} \mathds{1}_{g_{r,k} = 0}(P)\cdot\mathds{1}_{h_r = 0}(P)\\
        = {} & \frac{1}{q^2}\sum_{P\in\ff_q^{2r}}\sum_{a,b\in \ff_q}\psi(ag_{r,k}(P)+bh_r(P))\\
        = {} & \frac{1}{q^2}\sum_{a,b\in \ff_q}\sum_{P\in\ff_q^{2r}}\psi\left(\sum_{i=1}^k t_i(a,b) + \sum_{i = k+1}^r s_i (a,b) \right) \text{ where }P = (x_1,y_1,\ldots,x_r,y_r)\\
        = {} & \frac{1}{q^2}\sum_{a,b\in \ff_q}\sum_{x_i,y_i\in\ff_q} \prod_{i=1}^k \psi(t_i(a,b)) \prod_{i = k+1}^r \psi( s_i(a,b)) \text{ where } \psi(\alpha + \beta) = \psi(\alpha) \psi(\beta) \\
        = {} & \frac{1}{q^2}\sum_{a,b\in \ff_q} \prod_{i=1}^k \left(\sum_{x_i, y_i \in \ff_q} \psi(t_i(a,b))\right) \prod_{i = k+1}^r\left(\sum_{x_i,y_i\in\ff_q} \psi(s_i(a,b)) \right),
    \end{align*}
    where the last equality follows by summing independently over $x_i$ and $y_i$ for each $i$. Recall the definitions $S(a,b)$ and $T(a,b)$ in Section~\ref{sec3}. One obtains $$|V_{\ff_q}(g_{r,k},h_r)| = \frac{1}{q^2}\sum_{a,b\in \ff_q} T(a,b)^k S(a,b)^{r-k}.$$

    According to Lemma~\ref{lem: S(a,b)} and Lemma~\ref{lem: T(a,b)}, $$T(a,b)^k S(a,b)^{r-k} = \left\{\begin{array}{ll} q^{2r} & \text{if } a = b = 0, \\ q^r \psi(1)^k & \text{if } a \neq 0, b = 0, \\ (q(1 + N(a,b)))^r \psi(a^3 / b^2)^k & \text{if } b \neq 0. \end{array} \right.$$ As $q = 2^n$, one has $\psi(1) = (-1)^n$. Note that $N(a,b) \neq 2$ by Lemma~\ref{lem: mc{N}_k}. Moreover, it follows from Proposition~\ref{prop: determine Tr(a^3/b^2)} that for $b\neq 0$,
    $$\psi(a^3/b^2) = \begin{cases}
        (-1)^n &\text{ if }N(a,b) = 0\text{ or }3,\\
        (-1)^{n+1} &\text{ if }N(a,b) = 1.\\
    \end{cases}$$ 
    We obtain that $|V_{\ff_q}(g_{r,k},h_r)|$ is equal to \begin{align*} \frac{1}{q^2} \left( q^{2r} + (q-1) q^r (-1)^{nk} + q^r (-1)^{nk} \mc{N}_0 + (2q)^r (-1)^{(n+1)k} \mc{N}_1 + (4q)^r (-1)^{nk} \mc{N}_3 \right). \end{align*} Substituting the values of $\mc{N}_0$, $\mc{N}_1$, and $\mc{N}_3$ from Lemma~\ref{lem: mc{N}_k} completes the proof.
\end{proof}

\begin{cor}
\label{cor: characterize the equality of sizes}
Let the notation be as in Theorem~\ref{thm: sizes of rationals points}. For $k_1, k_2 \in \{0, 1, \ldots, r\}$, $$|V_{\ff_q}(g_{r,k_1}, h_r)| = |V_{\ff_q}(g_{r,k_2}, h_r)| \quad \Leftrightarrow \quad k_1 \equiv k_2 \pmod{2}.$$
\end{cor}

\begin{proof}
It is obvious that $k_1 \equiv k_2 \pmod{2}$ is a sufficient condition for the equality on the left-hand side. We explain the necessity. For convenience, let $a = \mc{I}_q(3) / q$, $b = 2^r \mc{I}_q(2)$, and $c = 4^r \binom{q}{3} / q$. Suppose the equality holds. If $n k_1 \equiv n k_2 \pmod{2}$ then $(-1)^{k_1} b = (-1)^{k_2} b$, which implies that both $k_1$ and $k_2$ have the same parity, since $b > 0$. If $n k_1 \not\equiv n k_2 \pmod{2}$, then one obtains $q-1 + a + c = 0$, which is impossible.
\end{proof}

\begin{remk}
\label{remk: bijection on rational points}
In view of Corollary~\ref{cor: characterize the equality of sizes}, it is natural to ask whether there is an explicit bijection between $V_{\ff_q}(g_{r,k_1}, h_r)$ and $V_{\ff_q}(g_{r,k_2}, h_r)$ when $k_1$ and $k_2$ have the same parity. This question reduces to constructing an explicit bijection between $V_{\ff_q}(g_{r,k}, h_r)$ and $V_{\ff_q}(g_{r,k+2}, h_r)$ for $0 \leq k \leq r-2$. We provide such a bijection below.

Consider the map
$$\phi_k:\ff_q^{2r}\to \ff_q^{2r};\quad \phi_k(x_1,y_1,\ldots,x_r,y_r) = (x_1',y_1',\ldots,x_r',y_r')$$
where
$$x_i' = \begin{cases}
    x_i &\text{ if } i\notin \{k+1,k+2\} ,\\
    x_i+z &\text{ if } i = k+1 \text{ or }k+2,\\
\end{cases}\qquad y_i' = \begin{cases}
    y_i &\text{ if } i\notin \{k+1, k+2\},\\
    y_i+z &\text{ if } i = k+1 \text{ or }k+2,
\end{cases}$$
and $$z = x_{k+1}+y_{k+1}+x_{k+2}+y_{k+2}.$$ Clearly, $\phi_k$ is an involution because $\ff_q$ has characteristic $2$. Moreover, one can check $g_{r,k} \circ \phi_k = g_{r, k+2}$ since $$x_{k+1}' y_{k+1}' + x_{k+2}' y_{k+2}' = x_{k+1}^2 + y_{k+1}^2 + x_{k+2}^2 + y_{k+2}^2 + x_{k+1} y_{k+1} + x_{k+2} y_{k+2},$$ and $h_r \circ \phi_k = h_r$ since \begin{align*} & x_{k+1}' y_{k+1}' (x_{k+1}' + y_{k+1}') + x_{k+2}' y_{k+2}' (x_{k+2}' + y_{k+2}') \\ = \ & x_{k+1} y_{k+1} (x_{k+1} + y_{k+1}) + x_{k+2} y_{k+2} (x_{k+2} + y_{k+2}).\end{align*} It follows that $\phi_k$ gives the desired bijection between $V_{\ff_q}(g_{r,k}, h_r)$ and $V_{\ff_q}(g_{r,k+2}, h_r)$.
\end{remk}

\section*{Acknowledgment}

The second author acknowledges support from National Science and Technology Council (NSTC) grants 114-2115-M-017-001-MY2 and 115-2115-M-017-007-MY3.



\bibliographystyle{alphaurl} 
\bibliography{Bib} 

\end{document}